\documentclass[]{amsart}
\pdfoutput=1

\usepackage{lscape}
\usepackage{subcaption}
\usepackage{scalerel}
\usepackage{mathrsfs}
\usepackage{hyphenat}
\usepackage{mlmodern}
\usepackage{mathpartir}

\usepackage[leftmargin=1em,rightmargin=1ex,vskip=1ex]{quoting}
\usepackage[utf8]{inputenc}
\usepackage{cmll}

\newcommand\scalemath[3]{\scalebox{#1}[#2]{\mbox{\ensuremath{\displaystyle #3}}}}

\newcommand{\leftarrowtip}{\ensuremath{\tikz\draw[line width=0.5pt,->] (10pt,0) -- (0,0);}}
\newcommand{\leftarrowtailnotip}{\ensuremath{\tikz\draw[line width=0.5pt,-<] (0,0) -- (10pt,0);}}

\newcommand{\unicodeStar}{\ensuremath{\star}}
\DeclareUnicodeCharacter{2605}{\unicodeStar}

\DeclareUnicodeCharacter{21D2}{\ensuremath{\Rightarrow}}
\DeclareUnicodeCharacter{2218}{\ensuremath{\circ}}
\DeclareUnicodeCharacter{2022}{\ensuremath{\bullet}}
\DeclareUnicodeCharacter{2219}{\ensuremath{\bullet}}
\DeclareUnicodeCharacter{2026}{\ensuremath{\dots}}
\DeclareUnicodeCharacter{2208}{\ensuremath{\in}}
\DeclareUnicodeCharacter{2192}{\ensuremath{\to}}
\DeclareUnicodeCharacter{2190}{\ensuremath{\leftarrowtip}}
\DeclareUnicodeCharacter{2919}{\ensuremath{\leftarrowtailnotip}}

\DeclareUnicodeCharacter{00D7}{\ensuremath{\times}}
\DeclareUnicodeCharacter{00B7}{\ensuremath{\cdot}}
\DeclareUnicodeCharacter{222B}{\ensuremath{\int}}
\DeclareUnicodeCharacter{22A4}{\ensuremath{\top}}
\DeclareUnicodeCharacter{22A5}{\ensuremath{\bot}}
\DeclareUnicodeCharacter{2264}{\ensuremath{\leq}}

\newcommand{\unicodecolon}{\ensuremath{\colon}}
\DeclareUnicodeCharacter{FE55}{\unicodecolon}
\newcommand{\unicodeleftpar}{\ensuremath{\left(}}
\DeclareUnicodeCharacter{27EE}{\unicodeleftpar}
\newcommand{\unicoderightpar}{\ensuremath{\right)}}
\DeclareUnicodeCharacter{27EF}{\unicoderightpar}
\DeclareUnicodeCharacter{2260}{\neq}
\DeclareUnicodeCharacter{22A9}{\Vdash}
\DeclareUnicodeCharacter{2237}{\proportion}
\DeclareUnicodeCharacter{2124}{\mathbb{Z}}
\DeclareUnicodeCharacter{27E8}{\langle}
\DeclareUnicodeCharacter{27E9}{\rangle}
\DeclareUnicodeCharacter{21A6}{\mapsto}
\DeclareUnicodeCharacter{22A2}{\vdash}
\DeclareUnicodeCharacter{2090}{\ensuremath{{}_a}}
\DeclareUnicodeCharacter{A71B}{{}^\uparrow}
\DeclareUnicodeCharacter{A71C}{{}^\downarrow}
\DeclareUnicodeCharacter{27E6}{\llbracket}
\DeclareUnicodeCharacter{27E7}{\rrbracket}

\newcommand{\unicoderightcircle}{\ensuremath{\RIGHTcircle}}
\DeclareUnicodeCharacter{25D1}{\unicoderightcircle}
\newcommand{\unicodeleftcircle}{\ensuremath{\LEFTcircle}}
\DeclareUnicodeCharacter{25D0}{\unicodeleftcircle}
\DeclareUnicodeCharacter{229B}{\circledast}

\newcommand{\unicodebbA}{\ensuremath{\mathbb{A}}}
\DeclareUnicodeCharacter{1D538}{\unicodebbA}
\newcommand{\unicodebbB}{\ensuremath{\mathbb{B}}}
\DeclareUnicodeCharacter{1D539}{\unicodebbB}
\newcommand{\unicodebbC}{\ensuremath{\mathbb{C}}}
\DeclareUnicodeCharacter{2102}{\unicodebbC}
\DeclareUnicodeCharacter{1D53B}{\ensuremath{\mathbb{D}}}
\DeclareUnicodeCharacter{2115}{\ensuremath{\mathbb{N}}}
\DeclareUnicodeCharacter{211D}{\ensuremath{\mathbb{R}}}
\DeclareUnicodeCharacter{1D543}{\ensuremath{\mathbb{L}}}
\newcommand\UnicodeBlackboardP{\ensuremath{\mathbf{P}}} \DeclareUnicodeCharacter{2119}{\UnicodeBlackboardP}
\DeclareUnicodeCharacter{211A}{\ensuremath{\mathbb{Q}}}
\DeclareUnicodeCharacter{1D544}{\ensuremath{\mathbb{M}}}
\DeclareUnicodeCharacter{1D54C}{\ensuremath{\mathbb{U}}}
\DeclareUnicodeCharacter{1D54D}{\ensuremath{\mathbf{V}}}
\DeclareUnicodeCharacter{1D54E}{\ensuremath{\mathbb{W}}}
\DeclareUnicodeCharacter{1D546}{\ensuremath{\mathbb{O}}}
\DeclareUnicodeCharacter{1D540}{\ensuremath{\mathbb{I}}}
\DeclareUnicodeCharacter{1D54A}{\ensuremath{\mathbb{S}}}
\DeclareUnicodeCharacter{1D53C}{\ensuremath{\mathbb{E}}}

\newcommand{\unicodecalS}{\ensuremath{\mathcal{S}}}
\newcommand{\unicodecalT}{\ensuremath{\mathcal{T}}}
\newcommand{\unicodecalC}{\ensuremath{\mathcal{C}}}
\newcommand{\unicodecalD}{\ensuremath{\mathcal{D}}}
\newcommand{\unicodecalX}{\ensuremath{\mathcal{X}}}
\newcommand{\unicodecalN}{\ensuremath{\mathcal{N}}}
\newcommand{\unicodecalE}{\ensuremath{\mathcal{E}}}
\DeclareUnicodeCharacter{1D4D4}{\unicodecalE}
\DeclareUnicodeCharacter{1D4D2}{\unicodecalC}
\DeclareUnicodeCharacter{1D4D3}{\unicodecalD}
\DeclareUnicodeCharacter{1D4DE}{\mathcal{O}}
\DeclareUnicodeCharacter{1D4AA}{\mathcal{O}}
\DeclareUnicodeCharacter{210B}{\mathcal{H}}
\DeclareUnicodeCharacter{1D4E2}{\unicodecalS}
\DeclareUnicodeCharacter{1D4E3}{\unicodecalT}
\DeclareUnicodeCharacter{1D4E7}{\unicodecalX}
\DeclareUnicodeCharacter{1D4D0}{\ensuremath{\mathcal{A}}}
\DeclareUnicodeCharacter{1D4D1}{\ensuremath{\mathcal{B}}}
\DeclareUnicodeCharacter{1D4D6}{\ensuremath{\mathcal{G}}}
\DeclareUnicodeCharacter{1D4D7}{\ensuremath{\mathcal{H}}}
\DeclareUnicodeCharacter{1D4DB}{\ensuremath{\mathcal{L}}}
\DeclareUnicodeCharacter{1D4DC}{\ensuremath{\mathcal{M}}}
\DeclareUnicodeCharacter{1D4DD}{\unicodecalN}
\DeclareUnicodeCharacter{1D4B1}{\ensuremath{\mathcal{V}}}
\DeclareUnicodeCharacter{1D4E5}{\ensuremath{\mathcal{V}}}
\DeclareUnicodeCharacter{1D4E6}{\ensuremath{\mathcal{W}}}
\DeclareUnicodeCharacter{1D4E4}{\ensuremath{\mathcal{U}}}

\DeclareUnicodeCharacter{03B1}{\alpha}
\DeclareUnicodeCharacter{03B2}{\beta}
\DeclareUnicodeCharacter{03BC}{\mu}
\DeclareUnicodeCharacter{03B4}{\delta}
\DeclareUnicodeCharacter{03B5}{\varepsilon}
\DeclareUnicodeCharacter{03B7}{\eta}
\DeclareUnicodeCharacter{03BB}{\lambda}
\DeclareUnicodeCharacter{03C1}{\rho}
\DeclareUnicodeCharacter{03C8}{\psi}
\DeclareUnicodeCharacter{03C4}{\tau}
\DeclareUnicodeCharacter{03A8}{\Psi}
\DeclareUnicodeCharacter{03C3}{\sigma}
\DeclareUnicodeCharacter{03C6}{\varphi}
\DeclareUnicodeCharacter{03A6}{\Phi}
\DeclareUnicodeCharacter{03A3}{\Sigma}
\DeclareUnicodeCharacter{03D5}{\phi}
\DeclareUnicodeCharacter{03B8}{\theta}
\DeclareUnicodeCharacter{03C0}{\ensuremath{\pi}}
\DeclareUnicodeCharacter{0393}{\Gamma}
\DeclareUnicodeCharacter{0394}{\Delta}
\DeclareUnicodeCharacter{03BA}{\kappa}
\DeclareUnicodeCharacter{03BD}{\nu}
\DeclareUnicodeCharacter{25A0}{\blacksquare}
\DeclareUnicodeCharacter{25AA}{\blacksquare}
\usepackage{amssymb}
\DeclareUnicodeCharacter{2205}{\varnothing}

\newcommand{\hirayo}{\scaleobj{0.9}{\text{\usefont{U}{min}{m}{n}\symbol{'210}}}}
\DeclareUnicodeCharacter{3088}{\hirayo}
\DeclareFontFamily{U}{min}{}
\DeclareFontShape{U}{min}{m}{n}{<-> udmj30}{}

\newcommand\UnicodeWhiteRightPointingSmallTriangle{\triangleright}
\DeclareUnicodeCharacter{25B9}{\mathbin{\UnicodeWhiteRightPointingSmallTriangle}}
\newcommand\UnicodeWhiteDownPointingSmallTriangle{\triangledown}
\DeclareUnicodeCharacter{25BF}{\mathbin{\UnicodeWhiteDownPointingSmallTriangle}}
\newcommand\UnicodeWhiteUpPointingSmallTriangle{\scalemath{1}{-1}{{}^{\triangledown}}}
\DeclareUnicodeCharacter{25B5}{\mathbin{\UnicodeWhiteUpPointingSmallTriangle}}

\DeclareUnicodeCharacter{2080}{\ensuremath{{}_0}}
\DeclareUnicodeCharacter{2081}{\ensuremath{{}_1}}
\DeclareUnicodeCharacter{2082}{\ensuremath{{}_2}}
\DeclareUnicodeCharacter{2083}{\ensuremath{{}_3}}

\DeclareUnicodeCharacter{1D62}{\ensuremath{{}_i}}
\DeclareUnicodeCharacter{2C7C}{\ensuremath{{}_j}}
\DeclareUnicodeCharacter{02B3}{\ensuremath{{}^r}}
\DeclareUnicodeCharacter{02E1}{\ensuremath{{}^\ell}}
\DeclareUnicodeCharacter{1D48}{\ensuremath{{}^d}}
\DeclareUnicodeCharacter{1D50}{\ensuremath{{}^m}}
\DeclareUnicodeCharacter{1D58}{\ensuremath{{}^u}}
\DeclareUnicodeCharacter{209A}{\ensuremath{{}_p}}
\DeclareUnicodeCharacter{2096}{\ensuremath{{}_k}}

\DeclareUnicodeCharacter{2245}{\ensuremath{\cong}}
\DeclareUnicodeCharacter{2286}{\subseteq}

\DeclareUnicodeCharacter{22C5}{\cdot}
\DeclareUnicodeCharacter{25C3}{\ensuremath{\triangleleft}}
\DeclareUnicodeCharacter{25B9}{\ensuremath{\triangleright}}

\DeclareUnicodeCharacter{2217}{\ast}
\DeclareUnicodeCharacter{039E}{\Xi}
\DeclareUnicodeCharacter{2295}{\oplus}
\DeclareUnicodeCharacter{2297}{\otimes}
\DeclareUnicodeCharacter{214B}{\parr}
\DeclareUnicodeCharacter{2298}{\oslash}
\DeclareUnicodeCharacter{25C0}{\mathbin{\blacktriangleleft}}
\DeclareUnicodeCharacter{25C1}{\mathbin{\vartriangleleft}}
\DeclareUnicodeCharacter{22B3}{\mathbin{\triangleright}}
\DeclareUnicodeCharacter{22B2}{\mathbin{\triangleleft}}
\DeclareUnicodeCharacter{FF5C}{\mid}
\DeclareUnicodeCharacter{227A}{\mathbin{\prec}}
\DeclareUnicodeCharacter{227B}{\mathbin{\succ}}
\DeclareUnicodeCharacter{22A3}{\mathbin{\dashv}}
\DeclareUnicodeCharacter{219D}{\ensuremath{\leadsto}}
\DeclareUnicodeCharacter{2191}{\ensuremath{\uparrow}}
\DeclareUnicodeCharacter{1361}{\colon}

\DeclareUnicodeCharacter{29D1}{\mathrel{\multimapdotbothB}}
\DeclareUnicodeCharacter{29D2}{\mathrel{\multimapdotbothA}}
\DeclareUnicodeCharacter{22C4}{\mathbin{\diamond}}
\DeclareUnicodeCharacter{226B}{\mathrel{\gg}}
\DeclareUnicodeCharacter{25A1}{\Box}
\DeclareUnicodeCharacter{266F}{\sharp}

\DeclareUnicodeCharacter{2099}{_n}
\DeclareUnicodeCharacter{2098}{_m}
\DeclareUnicodeCharacter{1D5AD}{\ensuremath{\mathsf{N}}}

\DeclareUnicodeCharacter{1D4DF}{\mathcal{P}}

\newcommand\mydots{\makebox[0.6em][c]{.\hfil.\hfil.}}
\DeclareUnicodeCharacter{2026}{\mydots}
\DeclareUnicodeCharacter{226B}{\gg}

\usepackage{stmaryrd}
\newcommand{\unicodeRelationalComposition}{\fatsemi}
\DeclareUnicodeCharacter{2A3E}{\unicodeRelationalComposition}
\usepackage{xcolor}

\definecolor{nordred}{HTML}{bf616a}
\definecolor{bordeaux}{HTML}{821529}
\definecolor{bluelink}{HTML}{003399}
\definecolor{nordred}{HTML}{bf616a}
\definecolor{nordblue}{HTML}{81a1c1}
\definecolor{norddarkblue}{HTML}{5e81ac}
\definecolor{nordgreen}{HTML}{a3be8c}
\definecolor{nordnight}{HTML}{4c566a}

\AtEndPreamble{
  \RequirePackage{hyperref}
  \RequirePackage{cleveref}
  \hypersetup{
    breaklinks = true,
    linktocpage,
    colorlinks = true,
    linkcolor = nordnight,
    citecolor = nordgreen,
    urlcolor = nordblue
  }
} %
\makeatletter
\newcommand{\nicelinktarget}[1]{\Hy@raisedlink{\hypertarget{#1}{}}}
\makeatother

\newcommand\defining[1]{\nicelinktarget{#1}{}}
\newcommand\Set{\hyperlink{linkSet}{\mathbf{Set}}}

\newcommand\comp{\ensuremath{\mathbin{⨾}}}

\newcommand\strictMonoidalCategory{\hyperlink{linkStrictMonoidalCategory}{strict monoidal category}}

\newcommand{\colorPre}[1]{{\color{nordred}{#1}}}
\newcommand{\colorMon}[1]{{\color{norddarkblue}{#1}}}

\newcommand\R{\colorPre{R}}

\newcommand\id{\mathrm{id}}

\usepackage[bibliography=common]{apxproof}
\renewcommand{\appendixsectionformat}[2]{Appendix for Section {#1}: {#2}}
\renewcommand{\appendixprelim}{\clearpage\appendix}

\usepackage[xcolor,no patch,hyperref,quotation,electronic]{knowledge}
\knowledgeconfigure{notion}
\knowledgestyle{notion}{color=nordnight}
\knowledgestyle{intro notion}{color=nordnight}

\knowledge{notion}
| effectful triple
| effectful triples
| Effectful triple
| Effectful triples

\knowledge{notion}
| double signature
| double signatures
| Double signature
| Double signatures

\knowledge{notion}
| forgetful double signature
| forgetful double signatures
| Forgetful double signature
| Forgetful double signatures

\knowledge{notion}
| double signature morphism
| double signature morphisms
| Double signature morphism
| Double signature morphisms
\newcommand{\doubleSig}{\kl[double signature morphism]{\mathsf{doubleSig}}}

\knowledge{notion}
| pinwheel double category
| pinwheel double categories
| Pinwheel double category
| Pinwheel double categories
| double category with pinwheels
| double categories with pinwheels
| Double category with pinwheels
| Double categories with pinwheels

\knowledge{notion}
| bigraph
| Bigraph
| bigraphs
| Bigraphs
| 2-graph
| 2-graphs
| 2-Graph
| 2-Graphs

\knowledge{notion}
| 2-graph morphism
| 2-graph morphisms
\newcommand{\Twograph}{\kl[2-graph morphism]{\mathsf{2\mbox{-}graph}}}
\newcommand{\Bigraph}{\kl[2-graph morphism]{\mathsf{2\mbox{-}graph}}}

\knowledge{notion}
| 2-graph category
| 2-graph categories
\newcommand{\Twocat}{\kl[2-category]{\mathsf{2\mbox{-}cat}}}

\knowledge{notion}
| double category
| double categories
| Double category
| Double categories

\knowledge{notion}
| weighted polygraph
| weighted polygraphs
| Weighted polygraph
| Weighted polygraphs
| timed polygraph
| timed polygraphs
| Timed polygraph
| Timed polygraphs

\knowledge{notion}
| morphism of timed polygraphs
| Morphism of timed polygraphs
| timed polygraph morphism
| Timed polygraph morphism
\newcommand{\tPolygraph}{\kl[morphism of timed polygraphs]{\mathsf{tPolygraph}}}

\knowledge{notion}
| tilted bicategory signature
| tilted bicategory signatures
| Tilted bicategory signature
| Tilted bicategory signatures
| tilted bigraph
| Tilted bigraph
| tilted bigraphs
| Tilted bigraphs
| tilted 2-graph
| Tilted 2-graph
| tilted 2-graphs
| Tilted 2-graphs
\newcommand{\tilt}{\kl[tilted bigraph]{\ensuremath{\mathsf{tilt}}}}

\knowledge{notion}
| path
| paths
| Path
| Paths
| composable path
| composable paths
\newcommand{\Path}{\kl[path]{\ensuremath{\mathsf{path}}}}

\knowledge{notion}
| signature
| signatures
| Signature
| Signatures

\knowledge{notion}
| congruence
| congruences
| Congruence
| Congruences

\knowledge{notion}
| symmetry
| symmetries
| Symmetry
| Symmetries

\knowledge{notion}
| profunctor
| profunctors
| Profunctor
| Profunctors

\knowledge{notion}
| left inclusion
| Left inclusion

\knowledge{notion}
| right inclusion
| Right inclusion

\knowledge{notion}
| left whiskering
| Left whiskering

\knowledge{notion}
| right whiskering
| Right whiskering

\knowledge{notion}
| interchange axiom
| Interchange axiom

\knowledge{notion}
| inclusion and whiskering
| Inclusion and whiskering
| whiskering interacts with inclusion
| Whiskering interacts with inclusion

\knowledge{notion}
| whiskering
| Whiskering

\knowledge{notion}
| rewiring preserves whiskering
| Rewiring preserves whiskering
| whiskering preserves rewiring
| Whiskering preserves rewiring

\knowledge{notion}
| rewiring preserves interchange
| Rewiring preserves interchange

\knowledge{notion}
| rewiring preserves composition
| Rewiring preserves composition

\knowledge{notion}
| rewiring
| Rewiring 
| definition of rewiring

\knowledge{rewiring is an action}[]{notion}

\knowledge{notion}
| composition
| Composition
| term composition
| Term composition

\knowledge{notion}
| term composition is associative
| Term composition is associative

\knowledge{notion}
| term composition is unital
| Term composition is unital

\knowledge{notion}
| substituting
| Substitution

\knowledge{notion}
| tensoring
| Tensoring

\knowledge{notion}
| arrow notation preterm
| Arrow notation preterm
| arrow notation preterms
| Arrow notation preterms

\knowledge{notion}
| arrow notation term
| Arrow notation term
| arrow notation terms
| Arrow notation terms

\knowledge{notion}
| observe-arrow notation
| Observe-arrow notation
| observe-arrow notation term
| Observe-arrow notation term
| observe-arrow notation terms
| Observe-arrow notation terms

\knowledge{notion}
| subdistribution
| Subdistribution
| subdistributions
| Subdistributions
| subdistributional

\knowledge{notion}
| rescaling
| Rescaling

\knowledge{notion}
| normalisation
| normalisations
| Normalisation
| Normalisations
| normalization
| normalizations
| Normalization
| Normalizations
| normalised
| Normalised

\knowledge{notion}
| observe-arrow notation copy-discard-compare category of terms

\knowledge{notion}
| category of arrow notation terms

\knowledge{notion}
| copy-discard-compare categories
| copy-discard-compare category
| Copy-discard-compare categories
| Copy-discard-compare category %
\theoremstyle{plain}
\newtheoremrep{theorem}{Theorem}
\newtheoremrep{proposition}[theorem]{Proposition}
\newtheoremrep{lemma}[theorem]{Lemma}
\newtheoremrep{corollary}[theorem]{Corollary}

\theoremstyle{definition}
\newtheorem{definition}[theorem]{Definition}

\theoremstyle{remark}
\newtheorem{remark}[theorem]{Remark}
\newtheorem{example}[theorem]{Example}

\allowdisplaybreaks

\title{Timing via Pinwheel Double Categories}
\author{Elena Di Lavore \and Mario Rom\'an}
\date{\today}

\makeatletter
\def\@copyrightspace{\relax}
\makeatother

\begin{document}

\begin{abstract}
  We discuss string diagrams for timed process theories—represented by
  duoidally-graded symmetric strict monoidal categories—built upon the string
  diagrams of pinwheel double categories.
\end{abstract}
\maketitle

\section{Timing process theories}

\begin{figure}[!ht]
  \includegraphics[width=0.5\textwidth]{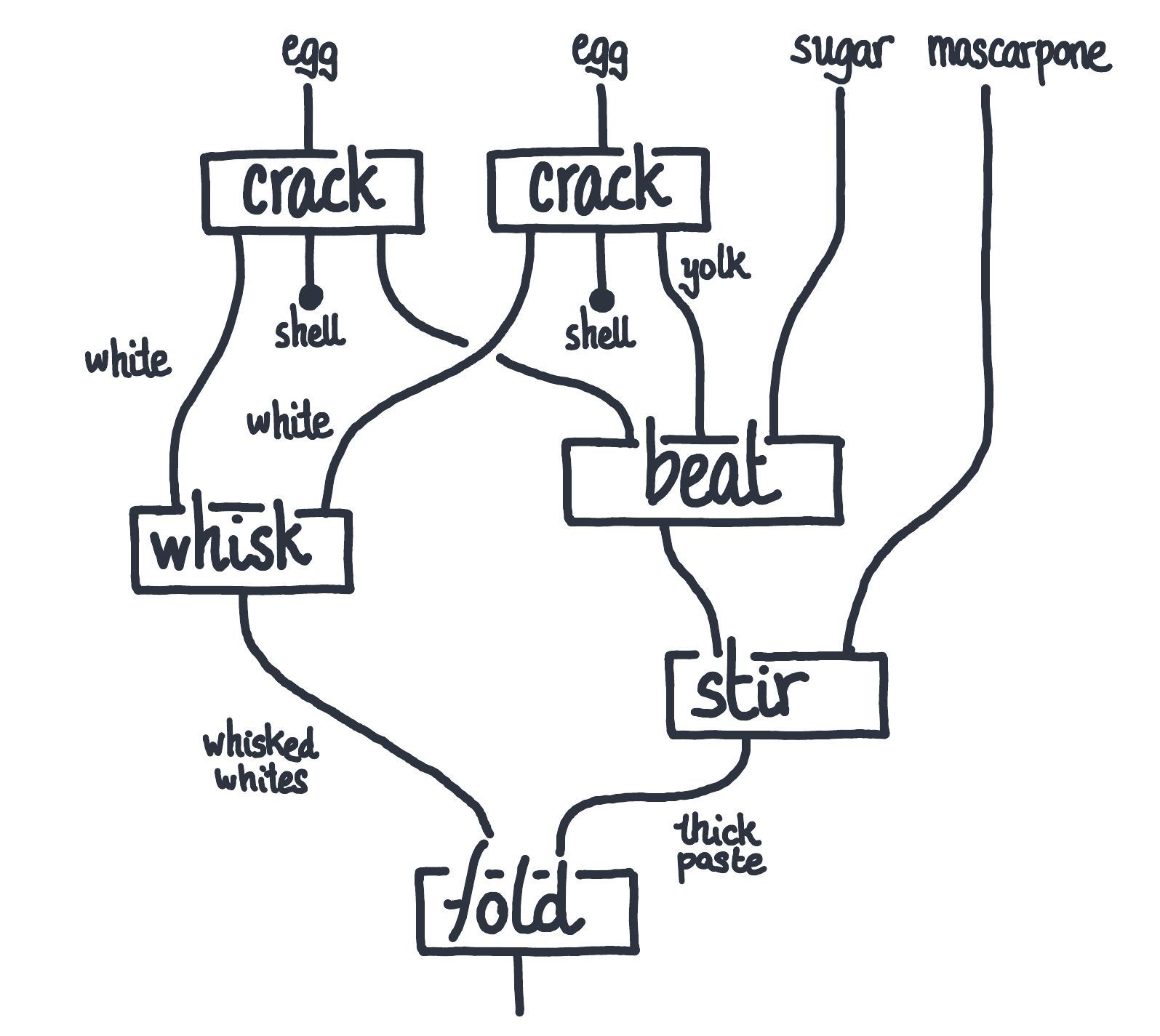}
  \caption{Process description of a preparation of \emph{crema di mascarpone}, adapted from Sobocinski}
  \label{fig:mascarpone}
\end{figure}

Process theories have an algebra in terms of symmetric monoidal categories; and
symmetric monoidal categories have an internal language in terms of string
diagrams, which we will use during this text. String diagrams are particularly
well-suited to process description (\Cref{fig:mascarpone}), and they appear
used, in different levels of formality, across science and engineering.

In symmetric monoidal categories, executing two independent processes in
parallel is the same as executing any of the two and then executing the other:
in fact, string diagrams cannot distinguish between the two. This is a feature,
for it reflects that the result must be the same in both cases. However, we may
be interested in aspects of the process beyond its result, and these may not be
preserved by naive string diagrammatic reasoning: for instance, parallelism
certainly matters if we are computing how much time a process takes.

\begin{figure}[!ht]
  \includegraphics[width=0.6\textwidth]{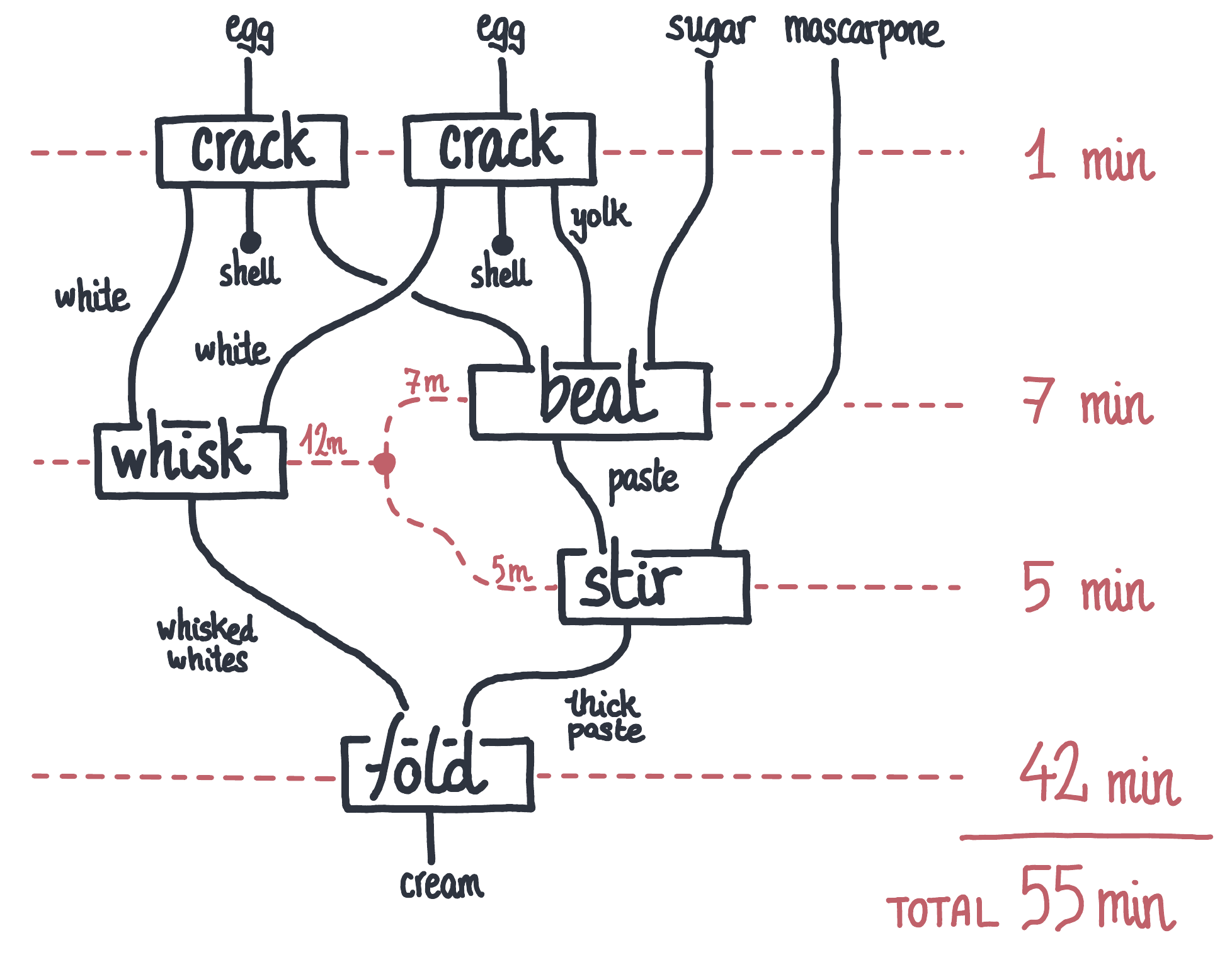}
  \caption{Timing of the ``crema di mascarpone'' process.}
  \label{fig:timing-mascarpone}
\end{figure}

This text proposes a string diagrammatic calculus for timing processes in monoidal
categories. It uses duoidal grading and string diagrams for pinwheel double categories.
While the underlying category theory needs some care, the resulting
string diagrams intuitively capture timing: in the vertical direction, we see how
resources get transformed, in the horizontal direction, we see how time flows.
Let us propose the following slogan.
\begin{quote}
  Process in sequence share \emph{resources}; process in parallel share \emph{time}.
\end{quote}

\subsection{Related work}
String diagrams for double categories \cite{dawson1995forbidden,myers2016string}
seem to still be missing a proof of initiality; only the single-object case has
been partially developed \cite{nester2021structure,nester2023concurrent}.
Instead, we take seriously the algebras of Delpeuch's monad on double signatures
\cite{delpeuch20:word} as a definition of pinwheel double categories. 
Duoidally graded string diagrams are missing, and
duoidal string diagrams are only developed for the normal case \cite{roman2024physicalduoidal}.

\section{Duoidal timing}

This section details the semantics of timed string diagrams in terms of
duoidally-graded monoidal categories—the monoidal counterpart of locally-graded
categories \cite{wood,mcDermott22}. The duoids we are interested in are
commutative and normal. They represent time durations, and they have two monoid
operations: an addition and a maximum. The main idea is that these operations
must satisfy the equation
$$
\mathsf{max}(a + c; b + d)
≤ 
\mathsf{max}(a; b) + \mathsf{max}(c; d).
$$
Interacting maximum and addition are a common motivation behind the development
of normal duoidal categories \cite{garner2016commutativity,shapiro2022duoidal}.

\begin{definition}[Duoid]
  A \emph{duoid}, $(A,≤, ⊕, 0, ↑, ⊥)$, is a poset $(A,≤)$ endowed with two monoid structures,
  $(A,⊕,0)$ and $(A,↑,⊥)$, that are monotone—meaning that $a ≤ a'$ and $b ≤ b'$ imply $a ⊕ a' ≤ b ⊕ b'$ and $a ↑ a' ≤ b ↑ b'$—and such that, additionally, the first laxly distributes over the second,
  $$(a ⊕ b) ↑ (c ⊕ d) ≤ (a ↑ c) ⊕ (b ↑ d),$$
  with $0 ↑ 0 ≤ 0$, and $⊥ ≤ ⊥ ⊕ ⊥$, and $⊥ ≤ 0$. The duoid is \emph{normal} when $⊥ = 0$. It is \emph{commutative} when $a ⊕ b = b ⊕ a$ and $a ↑ b = b ↑ a$.
\end{definition}

\begin{example}[Natural numbers]
  The natural numbers, with addition and the maximum, form a normal commutative
  duoid, $(ℕ,+,0,\mathsf{max},0)$. In the same sense that every monoidal
  category forms a duoidal category with the cartesian product when it exists,
  every monoidal poset forms a duoid with maximum when it exists. We are mostly
  concerned with these maximum-duoids.
\end{example}

\subsection{Duoidally-graded symmetric monoidal categories}
In a duoidally-graded symmetric monoidal category, every morphism $f ∈ ℂ(X;Y)$
has a grade, $a ∈ A$—we think of it as the time it takes to execute this
morphism—and we write $ℂ_a(X;Y)$ for the set of morphisms with grade $a ∈ A$. A
morphism can be regraded, forced to take more time: for each $f ∈ ℂ_a(X;Y)$
and $a ≤ b$, there exists $w_b(f) ∈ ℂ_b(X;Y)$.

\begin{toappendix}
\begin{definition}[Duoidally-locally graded monoidal category]
  For a \kl{duoidal category}, $(𝕍,⊗,I,⊲,N)$, a $𝕍$-graded (strict) monoidal
  category, $(ℂ,⊗,I)$, consists of (\emph{i}) a monoid of objects,
  $(ℂ_{obj},⊗,I)$; (\emph{ii}) a set of morphisms, $ℂ_A(X;Y)$, for each object
  $A ∈ 𝕍_{obj}$ and each pair of objects $X,Y ∈ ℂ_{obj}$; (\emph{iii}) a
  regrading vii 
  $$•[α] ፡ ℂ_A(X;Y) → ℂ_B(X;Y),\mbox{ for } α ∈ 𝕍(A;B),$$
  for each base morphism; (\emph{iv}) an identity morphism, $\id_X ∈ ℂ_N(X;X)$;
  (\emph{v}) a composition, 
  $$ (⨾) ፡ ℂ_A(X;Y) × ℂ_B(Y;Z) → ℂ_{A ⊲ B}(X;Z);$$
  and (\emph{vi}) a tensoring operation on morphisms,
  $$ (⊗) ፡ ℂ_A(X;Y) × ℂ_{A'}(X';Y') → ℂ_{A ⊗ A'}(X ⊗ X'; Y ⊗ Y');$$
  apart from (\emph{vii}) a special identity on the monoidal unit, $i_I ∈ ℂ_I(I;I)$.

  It must satisfy the following axioms, wherever well-typed.
  \begin{enumerate}
    \item $f ⨾ \id{} = f = \id{} ⨾ f$;
    \item $f ⨾ (g ⨾ h) = (f ⨾ g) ⨾ h$;
    \item $f ⊗ \id_I = f = \id_I ⊗ f$;
    \item $f ⊗ (g ⊗ h) = (f ⊗ g) ⊗ h$;
    \item $f[\id{}] = f$, and, in particular, $\id_X[\id_N] = \id_X$;
    \item $f[α ⨾ β] = f[α][β]$;
    \item $f[α] ⨾ g[β] = (f ⨾ g)[α ⊲ β]$;
    \item $f[α] ⊗ g[β] = (f ⊗ g)[α ⊗ β]$;
    \item $((f ⨾ g) ⊗ (f' ⨾ g'))[ψ] = (f ⊗ f') ⨾ (g ⊗ g')$;
    \item $i_I[ψ'] = \id_I$;
  \end{enumerate}  
  where we use the structure maps, $ψ ፡ (A ⊲ B) ⊗ (C ⊲ D) → (A ⊗ C) ⊲ (B ⊗ D)$ and $ψ' ፡ I → N$.
\end{definition}
\end{toappendix}

\begin{definition}[Duoid-graded monoidal category]
  For a \kl{duoid}, $(A,≤, ⊕, 0, ↑, ⊥)$, an $A$-graded (strict) monoidal
  category, $(ℂ,⊗,I)$, consists of (\emph{i}) a monoid of objects $(ℂ_{obj},⊗,I)$;
  (\emph{ii}) a set of morphisms, $ℂ_a(X;Y)$, for each grade $a ∈ A$ and each pair of
  objects $X,Y ∈ ℂ_{obj}$; (\emph{iii}) a regrading operation
  $$(•)[a ≤ b] ፡ ℂ_a(X;Y) → ℂ_b(X;Y)\mbox{, for each inequality }a ≤ b;$$
  (\emph{iv}) an identity morphism, $\id_X ∈ ℂ_0(X;X)$; (\emph{v}) a composition,
  $$(⨾) ፡ ℂ_a(X;Y) × ℂ_b(Y;Z) → ℂ_{a ⊕ b}(X;Z);$$
  and (\emph{vi}) a tensoring operation on morphisms
  $$
  (⊗) ፡ ℂ_a(X;Y) × ℂ_{a'}(X';Y') → ℂ_{a ↑ a'}(X ⊗ X'; Y ⊗ Y');
  $$
  apart from a special identity on the monoidal unit, $i_I ∈ ℂ_{⊥}(I;I)$.

  It must satisfy the following axioms, whenever well-typed.
  \begin{enumerate}
    \item $f ⨾ \id{} = f = \id{} ⨾ f$;
    \item $f ⨾ (g ⨾ h) = (f ⨾ g) ⨾ h$;
    \item $f ⊗ \id_I = f = \id_I ⊗ f$;
    \item $f ⊗ (g ⊗ h) = (f ⊗ g) ⊗ h$;
    \item $f[a ≤ a] = f$, and, in particular, $\id_X[⊥ ≤ ⊥] = \id_X$;
    \item $f[a ≤ c] = f[a ≤ b][b ≤ c]$;
    \item $f[a ≤ b] ⨾ g[c ≤ d] = (f ⨾ g)[a ⊕ c ≤ b ⊕ d]$;
    \item $f[a ≤ b] ⊗ g[c ≤ d] = (f ⊗ g)[a ↑ c ≤ b ↑ d]$;
    \item $((f ⨾ g) ⊗ (f' ⨾ g'))[(a ⊕ b) ↑ (c ⊕ d) ≤ (a ↑ c) ⊕ (b ↑ d)] = (f ⊗ f') ⨾ (g ⊗ g')$;
    \item $i_I[⊥ ≤ 0] = \id_I$;
  \end{enumerate}
\end{definition}

\begin{remark}
  We will be mostly concerned with symmetric monoidal categories graded over a
  normal and commutative maximum-duoid. Note how the interchange law does
  not hold automatically, $(f ⨾ h) ⊗ (g ⨾ k)$ takes equal or less time than $(f
  ⊗ g) ⨾ (h ⊗ k)$; but we can equate both by waiting during that difference of time.
\end{remark}

\section{Double categories}

\subsection{Double categories}
The usual definition of a \kl{double category} is that of a category internal to
the finitely complete category of small categories.

\begin{definition}[Double categories]
  A (strict) \emph{\intro{double category}} $𝔻$ consists of (\emph{1}) a set of
  objects, $𝔻_{obj}$; (\emph{2}) sets of arrows (or, \emph{horizontal arrows})
  $\smash{𝔻_{h}(X; Y)}$ for each pair of objects, $X,Y ∈ 𝔻_{obj}$;
  (\emph{3}) sets of arrows (or, \emph{vertical arrows})
  $\smash{𝔻_{v}(X; Y)}$ for each pair of objects, $X, Y ∈ 𝔻_{obj}$;
  and (\emph{4}) sets of cells, $𝔻(U;H;K;V)$, for each compatible quadruple of two
  vertical arrows $U ∈ 𝔻_{v}(X;Y)$ and $V ∈ 𝔻_{v}(Z;W)$, and two horizontal arrows
  $H ∈ 𝔻_{h}(X;Z)$ and $K ∈ 𝔻_{h}(Y;W)$.

  \kl{Double categories} have operations for vertical and horizontal composition of
  compatible 2-cells along their horizontal and vertical boundaries. These compositions
  are associative and unital and moreover respect the interchange law.
\end{definition}

\begin{definition}[Double signature]
  \AP A \emph{\intro{double signature}}, $Σ = (𝒪,ℋ,𝒱,𝓓)$, consists of (\emph{1}) a set
  of objects, $𝒪$; (\emph{2}) two graphs on these objects, given by sets $ℋ(X;Y)$ and
  $𝒱(X;Y)$ for each pair of objects, $X, Y ∈ 𝒪$; (\emph{3}) sets of generators,
  $$𝓓_{X,Y,Z,W}(u_1,...,u_n;h_1,...,h_p;k_1,...,k_q;v_1,...,v_m),$$
  for each four objects $X,Y,Z,W ∈ 𝒪$, and each quadruple of composable paths
  given by $\pmb{u} = u₁,...,u_n ∈ \Path(𝒱)(X;Y)$ and $\pmb{h} = h_1,...,h_p ∈
  \Path(ℋ)(Y;W)$, with $\pmb{v} = v_1,...,v_m ∈ \Path(Z;W)$, and $\pmb{k} =
  k_1,...,k_q ∈ \Path(ℋ)(X;Z)$.
\end{definition}

\begin{proposition}[Double signature morphism]
  \AP A \emph{\intro{double signature morphism}}, $α ፡ Σ → Σ'$, between
  two \kl{double signatures}, $Σ = (𝒪,ℋ,𝒱,𝓓)$ and $Σ' = (𝒪',ℋ',𝒱',𝓓')$,
  consists of 
  \begin{enumerate}
    \item a function on objects, $α_{obj} ፡ 𝒪 → 𝒪'$;
    \item two graph homomorphisms, with that same underlying function on
    objects, $α_H ፡ ℋ → ℋ'$ and $α_V ፡ 𝒱 → 𝒱'$, which extend to functors,
    $α_H^{∗}$ and $α_V^{∗}$; and
    \item a function on cells, preserving the boundary graphs,
    $$α ፡ 𝓓_{X,Y,Z,W}(\pmb{u};\pmb{h};\pmb{k};\pmb{v}) → 
    𝓓'_{α(X),α(Y),α(Z),α(W)}(α_V^{∗}(\pmb{u});α_H^{∗}(\pmb{h});α_H^{∗}(\pmb{k});α_V^{∗}(\pmb{v})).$$
  \end{enumerate}
  \kl{Double signatures} and \kl{double signature morphisms} form a category, $\doubleSig$.
\end{proposition}

\begin{proposition}[Forgetful double signatures]
  \AP The \intro{forgetful double signature} of a \kl{double category},
  $𝔻$, is a \kl{double signature} $Σ = (𝒪,ℋ,𝒱,𝓓)$ given by
  \begin{enumerate}
    \item the objects of the \kl{double category}, $𝒪 = 𝔻_{obj}$;
    \item the horizontal cells of the \kl{double category}, $ℋ(X;Y) = 𝔻_h(X;Y)$;
    \item the vertical cells of the \kl{double category}, $𝒱(X;Y) = 𝔻_v(X;Y)$;
    \item the 2-cells of the \kl{double category}, for each quadruple of composable paths
    given by $\pmb{u} = u₁,...,u_n ∈ \Path(𝒱)(X;Y)$ and $\pmb{h} = h_1,...,h_p ∈
    \Path(ℋ)(Y;W)$, with $\pmb{v} = v_1,...,v_m ∈ \Path(Z;W)$, and $\pmb{k} =
    k_1,...,k_q ∈ \Path(ℋ)(X;Z)$, we have
    $𝓓(\pmb{u}; \pmb{h}; \pmb{k}; \pmb{v}) = 
    𝔻(u_1 ⊗ ... ⊗ u_n; h_1 ⊗ ... ⊗ h_p; k_1 ⊗ ... ⊗ k_q; v_1 ⊗ ... ⊗ v_m).$
  \end{enumerate}
  The construction of the \kl{forgetful double signature} induces a functor from
  \kl{double categories} to \kl{double signatures} with strict double functors
  $\mathsf{forget} ፡ \mathsf{doubleCat} → \doubleSig$.
\end{proposition}

\subsection{Double categories from timed polygraphs}

Polygraphs are signatures for monoidal categories: they comprise some generating
objects and some basic processes from which all the others are built. If we
declare how much time these building blocks take, we can deduce how much every
composite process will take. This is the idea behind a \kl{timed polygraph}:
they are signatures for timed process theories.

\begin{definition}[Timed polygraph]
  \AP A \emph{\intro{timed polygraph}} $Ψ$ consists of a set of objects, $Ψ_obj$,
  and for each two lists of objects, $X_1,...,X_n ∈ Ψ_{obj}$ and $Y_1,...,Y_m ∈
  Ψ_{obj}$, and each natural number $t ∈ ℕ$, a set of generators,
  $Ψ_t(X_1,...,X_n; Y_1,...,Y_m)$.
\end{definition}

\begin{definition}[Timed polygraph morphism]
  \AP A \emph{\intro{timed polygraph morphism}}, $α ፡ Ψ → Ξ$, consists of a function
  between the object sets, $α_{obj} ፡ Ψ_{obj} → Ξ_{obj}$, and a function between generators
  that preserves time, inputs and outputs,
  $$α ፡ Ψ_{t}(X_1,...,X_n; Y_1,...,Y_m) → Ξ_{t}(α(X_1),...,α(X_n); α(Y_1),...,α(Y_m)).$$
\end{definition}

\begin{figure}[!ht]
  \includegraphics[width=0.7\textwidth]{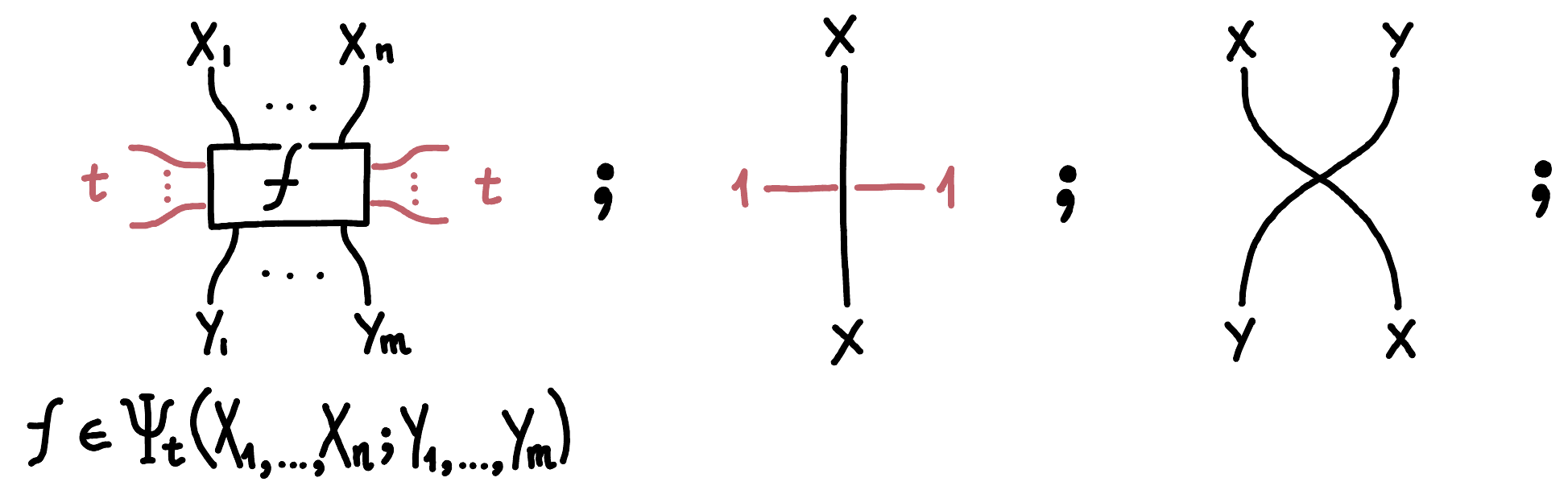}
  \caption{Double signature from a timed polygraph.}
  \label{fig:timing-signature}
\end{figure}

\begin{proposition}[Double signature from a timed polygraph]
  Every \kl{timed polygraph} $Ψ$ induces a \kl{double signature},
  $\mathsf{draw}(Ψ) = (𝒪,ℋ,𝒱,𝓓)$, containing
  \begin{enumerate}
    \item a single object, $𝒪 = \{ ∗ \}$;
    \item a single horizontal generator, $ℋ(∗,∗) = \{ 1 \}$;
    \item a vertical generator for each object, $𝒱(∗,∗) = Ψ_{obj}$;
    \item and three kinds of cells: a cell for each timed generator, 
    $$𝓓(X_1,...,X_n;t;t;Y_1,...,Y_m) = Ψ_t(X_1,...,X_n;Y_1,...,Y_m),$$
    and, additionally, ``waiting'' cell, $𝓓(X;1;1;X) = Ψ_1(X;X) + \{\mathrm{wait}\}$, and a ``braiding'' cell
    $𝓓(X,Y;0;0;Y,X) = Ψ_0(X,Y;Y,X) + \{σ\}$ (see \Cref{fig:timing-signature}).
  \end{enumerate}
  This assignment extends to a functor $\mathsf{draw} ፡ \tPolygraph → \doubleSig$.
\end{proposition}

\begin{figure}[!ht]
  \includegraphics[width=0.75\textwidth]{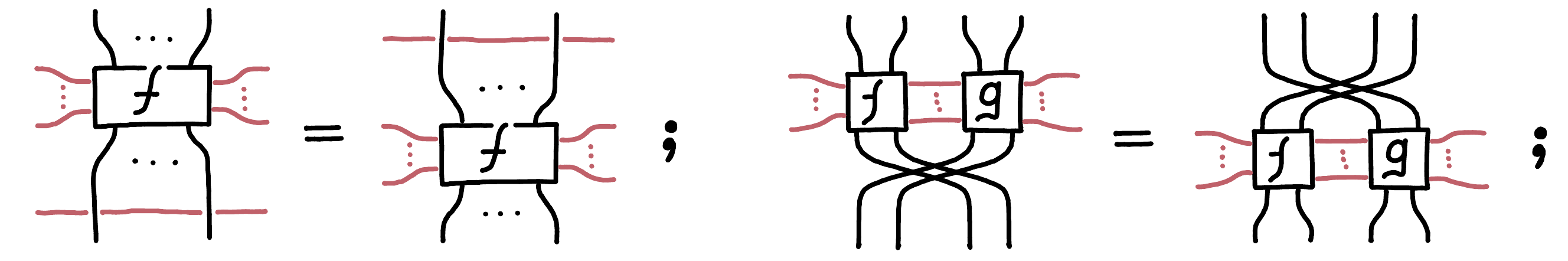}
  \caption{Equations for the double category from a timed polygraph.}
  \label{fig:timing-equations}
\end{figure}

\begin{proposition}[Double category from a timed polygraph]
  Every \kl{timed polygraph}, $Ψ$ induces a \kl{double category} freely
  generated by its associated \kl{double signature}, $\mathsf{draw}(Ψ)$, and
  quotiented by the equations in \Cref{fig:timing-equations}.
\end{proposition}

\subsection{Towards pinwheel double categories}
String diagrams for double categories face a problem: not every plane
arrangement of composable cells can be explained in terms of the algebra of
double categories \cite{dawson1995forbidden}. Those that can be explained were described by Dawson as
``neat tilings'' and characterized as those that do not contain a full
\emph{pinwheel} (\Cref{fig:pinwheels}): the only obstruction to composition in double categories
are pinwheels.

\begin{example}
  Consider a quintuple of timed processes: $f$ and $h$ take two units of time,
  while $g$, $k$, and $a$ take a single unit. The resources of $a$ and $h$
  depend on $g$, while the resouces of $k$ depend on $f$ and $a$. 
  We could: %
  \emph{(1)} execute $f$ and $g$ in parallel; %
  \emph{(2)} when $g$ finishes, execute $a$ and $h$ in parallel; %
  \emph{(3)} when both $f$ and $a$ finish, execute $k$; %
  and \emph{(4)} both $h$ and $k$ will finish at the same time, after 4 units of time.
  There is no physical limitation for this arrangement, but it is disallowed by
  the algebra of double categories (\Cref{fig:pinwheels}, left).

\begin{figure}[!ht]
  \includegraphics[width=0.5\textwidth]{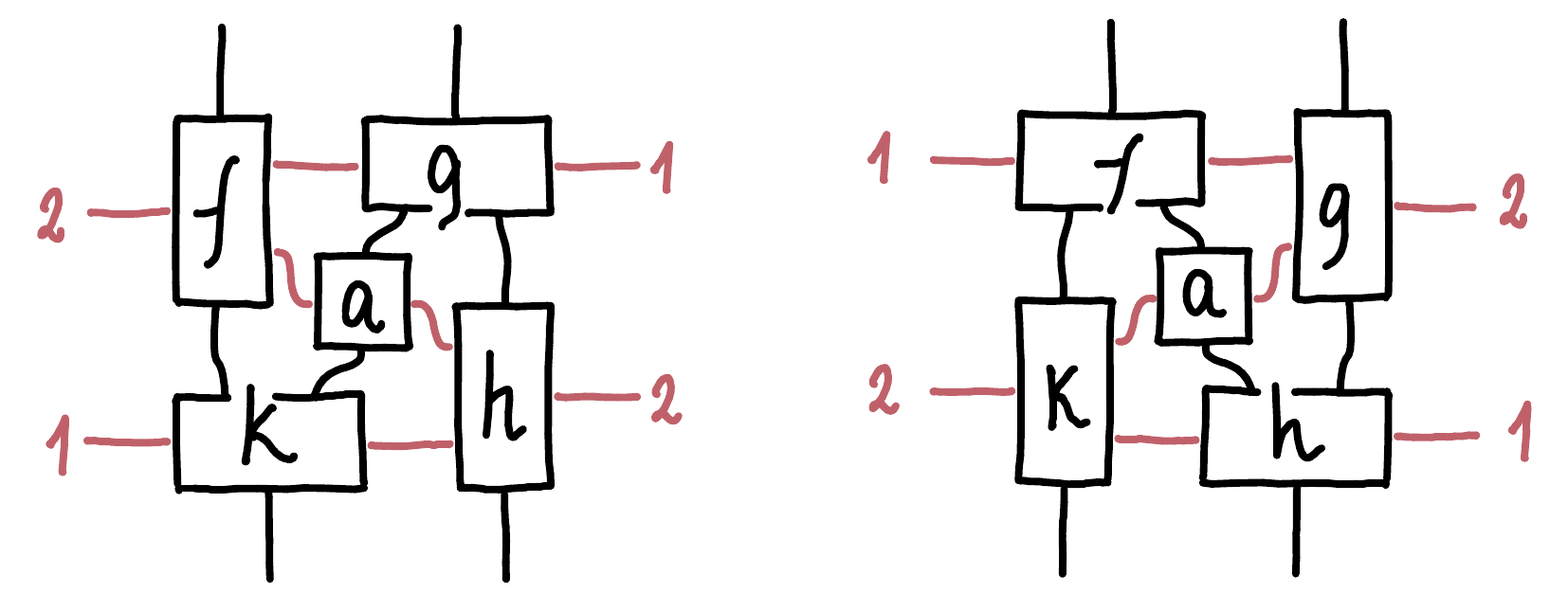}
  \caption{Two valid timings that form pinwheels.}
  \label{fig:pinwheels}
\end{figure}
\end{example}

However, this means double categories cannot be the algebraic structure that
encodes all possible combinations of timing processes. The pinwheel represents a
valid timing, and so we should account for it. Instead, we need more algebraic
structure: we need double categories \emph{with pinwheels}. Recent work by
Delpeuch will enable this route \cite{delpeuch20:word}.

\section{Double categories with pinwheels}%
\label{sec:doubleWithPinwheels}%

\begin{figure}[!ht]
  \includegraphics[width=0.5\textwidth]{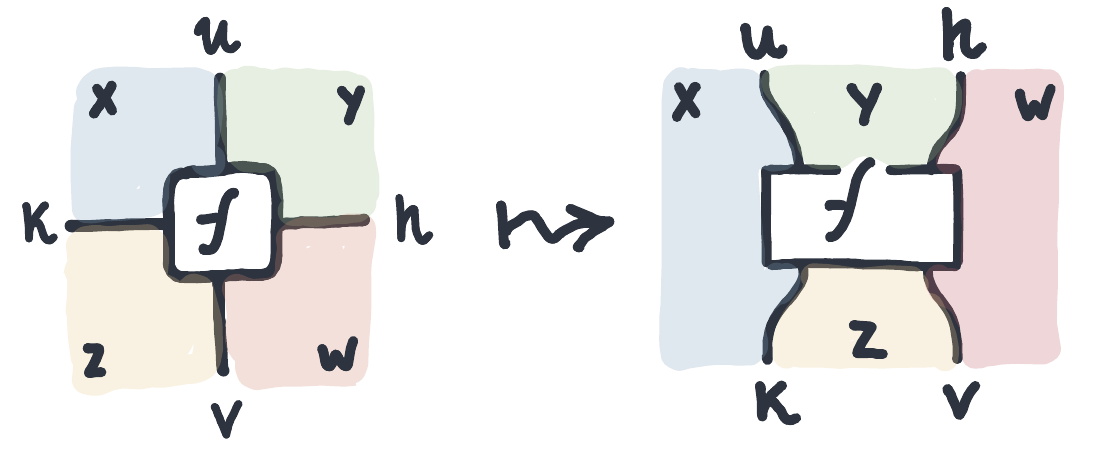}
  \caption{Double cells can be tilted to 2-cells.}
  \label{fig:tilted-double-cell}
\end{figure}

All uses the double categories that do not rely on the exclusion of pinwheels
can be recasted for 2-categories, gaining the  extra expressivity for pinwheels
in the process \cite{delpeuch20:word}. The main idea is that any 2-cell in a
double category can be ``tilted'' and interpreted as a 2-cell in a corresponding
2-category that has, as generating objects, both the vertical and horizontal
cells of the double category (\Cref{fig:tilted-double-cell}).

In formal terms, each \kl{double signature} can be translated into a
\kl{2-graph} that generates the free double category with pinwheels.

\begin{definition}[Tilted 2-graph]
  The \emph{\intro{tilted 2-graph}} of a \kl{double signature} $Σ =
  (𝒪,ℋ,𝒱,𝓓)$ is the \kl{2-graph}, $\tilt(Σ)$, defined as having
  \begin{itemize}
    \item the same 0-cells, $\tilt(Σ)_{obj} = 𝒪$;
    \item both vertical and horizontal 1-cells, $\tilt(Σ)(X;Y) = 𝒱(X;Y) +
    ℋ(X;Y)$;
    \item the same 2-cells between $X ∈ 𝒪$ and $W ∈ 𝒪$ when paths are composed
    precisely of $\pmb{u} = u_1,...u_n ∈ \Path(𝒱)(X;Y)$ and $\pmb{h} =
    h_1,...,h_p ∈ \Path(ℋ)(Y;W)$, with $\pmb{v} = v_1,...,v_m ∈ \Path(𝒱)(Z;W)$
    and $\pmb{k} = k_1,...,k_q ∈ \Path(ℋ)(X;Z)$, for any two $Y, Z ∈ 𝒪$; that
    is,
    $$
    \qquad\qquad \tilt(Σ)_{X,W}(u_1,...,u_n,h_1,...,h_p; k_1,...,k_q,v_1,...,v_n) = 
    𝓓_{X,Y,Z,W}(\pmb{u};\pmb{h};\pmb{k};\pmb{v});
    $$
    and empty otherwise, with $\tilt(Σ)(x_1,...,x_n;y_1,...,y_m) = ∅$ for any
    other two given paths $\pmb{x} = x_1,...,x_n ∈ \Path(\tilt(Σ))(X;W)$ and $\pmb{y}
    = y_1,...,y_m ∈ \Path(\tilt(Σ))(X;W)$.
  \end{itemize}
\end{definition}

From this \kl{2-graph}, we can construct its free \kl{2-category}, 
$\mathsf{string}(\tilt(Σ))$. String diagrams for this category contain the cells
of the pinwheel double category we want to construct, but they also contain many
diagrams that cannot be interpreted as double cells: we use those that do to
construct the tentative \emph{free \kl{pinwheel double category}} over a
\kl{double signature}.

\begin{figure}[!ht]
  \includegraphics[width=0.5\textwidth]{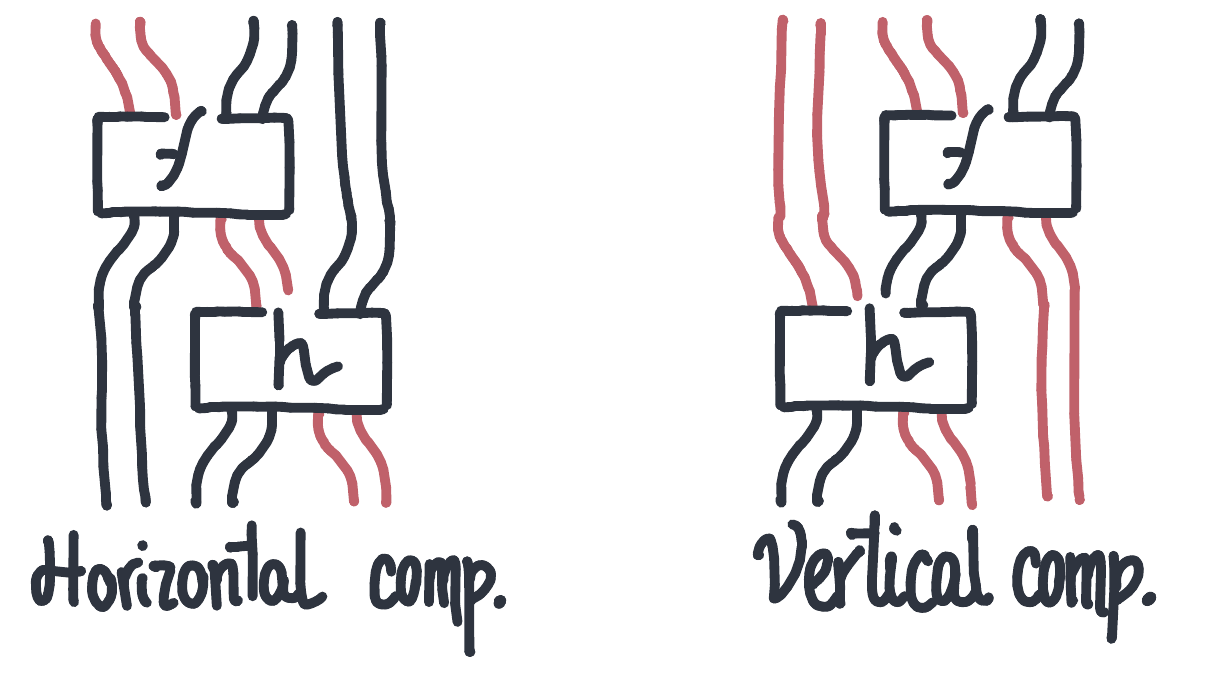}
  \caption{Horizontal and vertical compositions in the tilted 2-category.}
  \label{fig:tilted-compositions}
\end{figure}

\begin{definition}[Pinwheel double category]
  Any \kl{double signature}, $Σ = (𝒪,ℋ,𝒱,𝓓)$, generates a strict \kl{double category},
  $\mathsf{pinwheel}(Σ)$, having
  \begin{enumerate}
    \item objects from its signature $\mathsf{pinwheel}(Σ)_{obj} = 𝒪$;
    \item horizontal paths, composing as such, $\mathsf{pinwheel}(Σ)_{h} = \Path(ℋ)$;
    \item vertical paths, composing as such, $\mathsf{pinwheel}(Σ)_{v} = \Path(𝒱)$;
    \item cells the string diagrams constructed from the \kl{tilted 2-graph};
    that is, for $\pmb{u} = u_1,...u_n ∈ \Path(𝒱)(X;Y)$ and $\pmb{h} =
    h_1,...,h_p ∈ \Path(ℋ)(Y;W)$, with $\pmb{v} = v_1,...,v_m ∈ \Path(𝒱)(Z;W)$
    and $\pmb{k} = k_1,...,k_q ∈ \Path(ℋ)(X;Z)$, we have that the set of cells
    $\mathsf{pinwheel}(Σ)(\pmb{u};\pmb{h};\pmb{k};\pmb{v})$ is defined by
    $$\mathsf{string}(\tilt(Σ))(u_1 ⊗ ... ⊗ u_n ⊗ h_1 ⊗ ... ⊗ h_p; k_1 ⊗ ... ⊗ k_q ⊗ v_1 ⊗ ... ⊗ v_m).$$
  \end{enumerate}
  Horizontal and vertical composition for this double category are given by compositions of 2-categorical
  string diagrams, following \Cref{fig:tilted-compositions}.
\end{definition}

\begin{proposition}[Pinwheel monad]
  \AP Pinwheel forms a monad
  $$\mathsf{Pinwheel} = (\mathsf{forget} ⨾ \mathsf{pinwheel}) ፡ \doubleSig → \doubleSig.$$
\end{proposition}

\begin{definition}[Pinwheel double category]
  \AP \intro{Pinwheel double categories} are the algebras of the pinwheel monad on
  \kl{double signatures}.
\end{definition}

\begin{toappendix}
\section{String Diagrams of Bicategories}
\label{sec:strings-bicategories}

\begin{definition}
  \defining{linktwocategory}{}
  A strict \emph{2-category} $𝔹$ consists of a collection of \emph{objects}, or
  0-cells, $𝔹_{obj}$, and a category of \emph{morphisms} or 1-cells between any
  two objects, $𝔹(A;B)$. A strict \emph{2-category} is endowed with operations
  for the parallel composition of 1-cells,
  \begin{align*}
    (\comp) &፡ 𝔹(A; B) × 𝔹(B; C) → 𝔹(A; C), \\
    (I_A) &፡ 𝔹(A;A),
  \end{align*}
  that are associative and unital both on objects and morphisms, meaning that
  $(X \comp Y) \comp Z = X \comp (Y \comp Z)$, and $I_A \comp X = X = X \comp I_B$.
  Bicategories must satisfy the following axioms, making parallel composition a functor:
  \begin{enumerate}
    \item parallel composition is unital, $f \comp \id = f$, and $\id \comp f = f$;
    \item parallel composition is associative, $f \comp (g \comp h) = (f \comp g) \comp h$;
    \item compositions are unital, $\id \comp \id = \id$;
    \item compositions interchange, $(f ⨾ g) \comp (f' ⨾ g') = (f \comp f') ⨾ (g \comp g')$.
  \end{enumerate}
\end{definition}
\begin{remark}
  A single-object strict 2-category is exactly a \strictMonoidalCategory{}.
\end{remark}

\subsection{String diagrams of 2-categories}
Let us assume any construction of bicategorical string diagrams. There are
multiple options when constructing bicategorical string diagrams: we could
propose a combinatorial description or a topological one, we could decide that
0-cells appear as part of the structure, or simply as the property of being
well-typed. These constructions are all left adjoints to the same forgetful functor
from 2-categories to \emph{\kl{2-graphs}}.

\begin{definition}[2-graph]
  \AP A \emph{\intro{2-graph}}, or \emph{bigraph}, $𝓑$ is given by a set of objects, $𝓑_{obj}$; a set of arrows between any two objects, $𝓑(A; B)$; and a set of 2-arrows between any two paths of arrows, $𝓑(X₀,…,Xₙ; Y₀,…,Yₘ)$.
\end{definition}

\begin{definition}[2-graph morphism]
  \AP A \emph{\intro{2-graph morphism}}, $f ፡ 𝓐 → 𝓑$, is a function between
  their object sets, $f_o ፡ 𝓐_{obj} → 𝓑_{obj}$; a family of functions between
  their corresponding arrow sets, $f ፡ 𝓐(A;B) → 𝓑(f(A),f(B))$; and a family of
  functions between their corresponding 2-arrow sets,
  $$f ፡ 𝓐(X_0,…,X_n; Y_0,…,Y_m) → 𝓑(f(X_0),…,f(X_n); f(Y_0),…,f(Y_m)).$$
  \kl{2-Graphs} with \kl{2-graph morphisms} form a category, $\Bigraph$.
\end{definition}

\begin{definition}[Forgetful 2-graph]
  The forgetful \kl{2-graph}, $\mathsf{forget}(𝔸)$, of a \kl{2-category}, $𝔸$,
  has the same set of objects as the \kl{2-category}, $\mathsf{forget}(𝔸)_{obj}
  = 𝔸_{obj}$; and morphism sets given by all the morphisms of the \kl{2-category},
  $$
  \mathsf{forget}(𝔸)(X_1,...,X_n;Y_1,...,Y_n) =
  𝔸(X_1 ⊗ ... ⊗ X_n; Y_1 ⊗ ... ⊗ Y_n).
  $$
\end{definition}

\begin{theorem}[Free bicategories]
  The functor $\mathsf{forget} ፡ \Twocat → \Twograph$ is a right adjoint.
\end{theorem}

That is, there is a right adjoint from \kl{2-categories} to \kl{2-graphs}. A
left adjoint can be constructed by bicategorical string diagrams over the
bigraph; let us write $\mathsf{string} ፡ \Twograph → \Twocat$ for this left
adjoint. However, let us still not assume any particular construction for this
left adjoint.
\end{toappendix}

\section{Conclusions}

We have seen how to construct timed string diagrams, like that in
\Cref{fig:timing-mascarpone}, using \kl{pinwheel double categories}.
In this example, the \kl{timed polygraph} is given by the following
generators.

\begin{figure}[!ht]
  \includegraphics[width=0.6\textwidth]{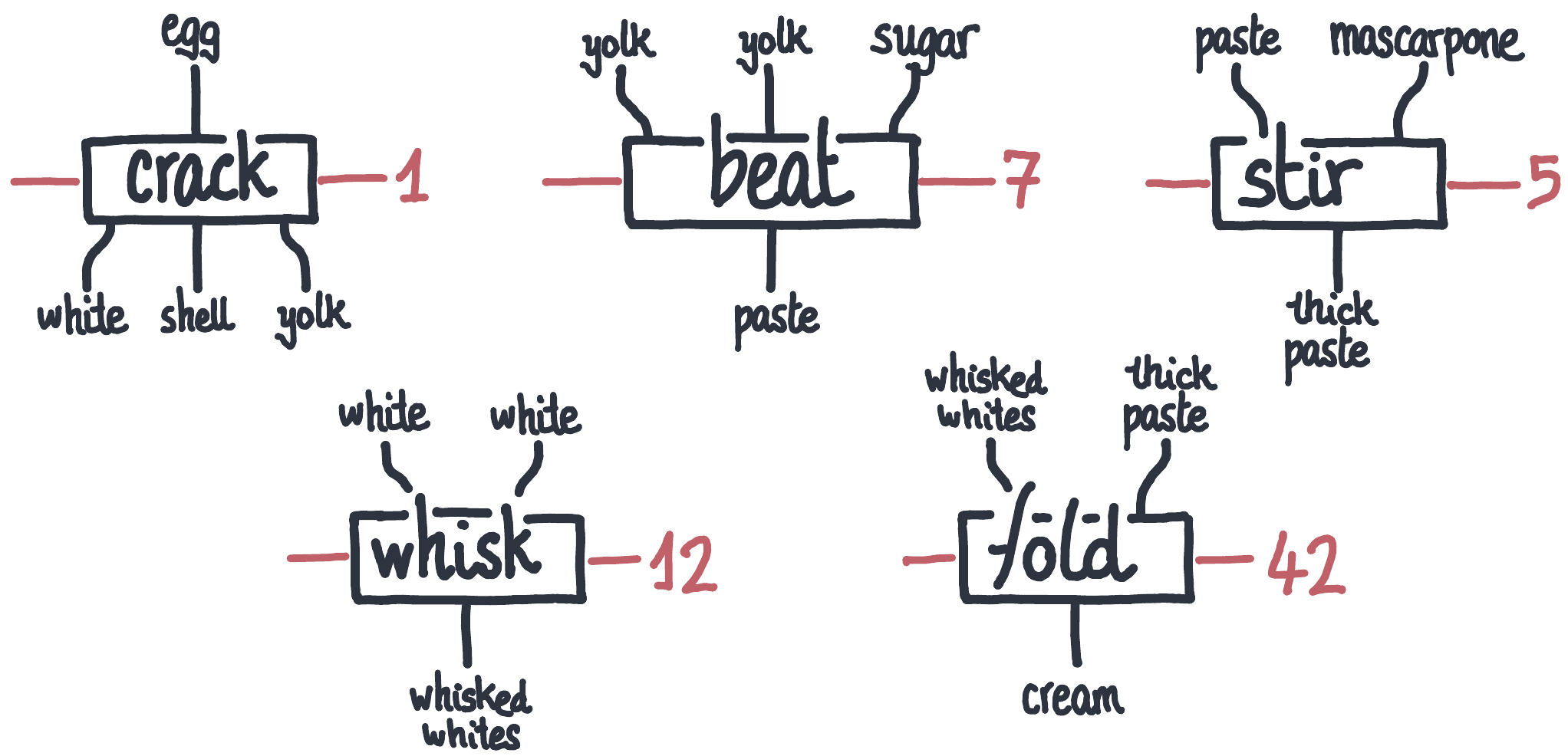}
  \caption{Timed polygraph for crema di mascarpone.}
  \label{fig:timed-signature}
\end{figure}

\subsection{Further work}
More details for the constructions of this paper may be provided in later
versions. Monoidal width, developed by Di Lavore and Sobocinski, may arise from
a slight generalization of the construction presented here
\cite{diLavoreSobocinski23:monoidalWidth}. It may be possible to connect this way
of timing processes to time-graded coalgebras of monads, which the authors have recently
proposed for continuous dynamical systems \cite{dilavore:graded}. 

\subsection*{Funding}
Mario Román was supported by the Air Force Office of Scientific Research under
award number FA9550-21-1-0038. Elena Di Lavore and Mario Román were supported by
the Advanced Research + Invention Agency (ARIA) Safeguarded AI Programme.

\bibliographystyle{alpha}
\bibliography{main.bib}

\end{document}